\documentclass[a4paper,pdftex,reqno,10pt]{amsart}

\allowdisplaybreaks[1]

\usepackage{geometry}
\usepackage{graphicx}
\usepackage{enumerate}
\usepackage{courier}
\usepackage{times}
\usepackage{amssymb}
\usepackage{amsmath}

\usepackage{hyperref}

\theoremstyle{plain}
\newtheorem{Theorem}{Theorem}[section]
\newtheorem{Proposition}[Theorem]{Proposition}
\newtheorem{Conjecture}[Theorem]{Conjecture}
\newtheorem{Corollary}[Theorem]{Corollary}
\newtheorem{Lemma}[Theorem]{Lemma}

\newenvironment{Proof}
{\begin{trivlist}\item[]{{\sc Proof.}}}{\hfill{$\square$}\noindent\end{trivlist}}

\newcommand{\gaussm}[3]{\genfrac{[}{]}{0pt}{}{#1}{#2}_{#3}}

\newcommand{\F}[2]{\mathbb{F}_{#2}^{#1}}

\theoremstyle{definition}
\newtheorem{Definition}[Theorem]{Definition}

\theoremstyle{remark}

\begin{document}

%%%%%%%%%%%%%%%%%%%%%%%%%%%%%%%%%%%%%%%%%%%%%%%%%%%%%%%%%%%%%%%%%%%%%%%%%%%%%%
%% Title:
%%%%%%%%%%%%%%%%%%%%%%%%%%%%%%%%%%%%%%%%%%%%%%%%%%%%%%%%%%%%%%%%%%%%%%%%%%%%%%

\title{Generalized vector space partitions}

\author[Daniel Heinlein]{Daniel Heinlein}
\author[Thomas Honold]{Thomas Honold}
\author[Michael Kiermaier]{Michael Kiermaier}
\author[Sascha Kurz]{Sascha Kurz}
\address{Department of Mathematics, University of Bayreuth, 95440 Bayreuth, Germany.}
\email{daniel.heinlein@uni-bayreuth.de}
\email{michael.kiermaier@uni-bayreuth.de}
\email{sascha.kurz@uni-bayreuth.de}
\address{ZJU-UIUC Institute, Zhejiang University, 314400 Haining, China.}
\email{honold@zju.edu.cn}
\thanks{$^\star$ The work of the authors was partially supported by the grants KU 2430/3-1, WA 1669/9-1 -- 
Integer Linear Programming Models for Subspace Codes and Finite Geometry from the German Research Foundation.}

\date{}

\maketitle

\begin{abstract}
  A \emph{vector space partition} $\mathcal{P}$ in $\F{v}{q}$ is a set of subspaces such that every 
  $1$-dimensional subspace of $\F{v}{q}$ is contained in exactly one element of $\mathcal{P}$. Replacing {\lq\lq}$1$-dimensional{\rq\rq} by
  {\lq\lq}$t$-dimensional{\rq\rq}, we generalize this notion to vector space $t$-partitions and study their properties.
  There is a close connection to subspace codes and some problems are even interesting and unsolved for the set case $q=1$.  
  
  \medskip
  
  \noindent
  \textbf{Keywords:} Galois geometry, partial spreads, constant-dimension codes, subspace codes, $q$-analogs, pairwise balanced designs,  
  and vector space partitions\\
  \textbf{MSC:} 51E23; 05B40 
  %%\textbf{MSC:} 51E23; 05B15, 05B40, 11T71, 94B25 
  %% 51E23 Spreads and packing problems
  %% 94B25 Combinatorial codes
  %% 05B15 Orthogonal arrays, Latin squares, Room squares
  %% 05B40 Packing and covering
  %% 11T71 Algebraic coding theory; cryptography 
\end{abstract}

\section{Introduction}
\label{sec_introduction}

\noindent
A vector space partition of $\F{v}{q}$ consists of subspaces such that every
$1$-dimensional subspace is covered exactly once. As a natural 
extension we consider sets of subspaces such that every $t$-dimensional subspace is
covered exactly once. For $t\ge 2$ this leads
to questions that even remain unsolved for the set case $q=1$. For $q\ge 2$ there is
a close relation to constant-dimension or 
mixed-dimension subspace codes with respect to the subspace distance. Having such a
code at hand, intersecting every codeword with 
a hyperplane gives an object that cannot be described as a mixed-dimension subspace
code in terms of the minimum subspace distance 
directly. However, our generalization of a vector space partition captures this
situation and yields non-trivial upper bounds 
for constant-dimension subspace codes.  

More precisely, let $q>1$ be a prime power and $v$ a positive integer. A
\emph{vector space partition} $\mathcal{P}$ of $\F{v}{q}$ is a set of 
subspaces with the property that every non-zero vector is contained in a unique
member of $\mathcal{P}$. If $\mathcal{P}$  contains 
$m_d$ subspaces of dimension $d$, then $\mathcal{P}$ is of \emph{type} $v^{m_v}\dots
1^{m_1}$. We may leave out %%some of the 
dimensions with $m_d=0$. Subspaces of dimension~$1$ in $\mathcal{P}$ are called
\emph{holes}.
%If there is at least one non-hole, then $\mathcal{P}\neq\{\F{v}{q}\}$ is called non-trivial. 
The vector space partition consisting only of holes and the vector space partition
$\{\F{v}{q}\}$ are called \emph{trivial}.

Here, we give a natural generalization of this notion. For a positive integer $t$,
we call a set $\mathcal{P}$ of subspaces of $\F{v}{q}$ a \emph{vector space
$t$-partition}, if all elements of $\mathcal{P}$ are of dimension at least $t$ and
every $t$-dimensional subspace is contained in a unique element of $\mathcal{P}$.
Ordinary vector space partitions are precisely the vector space $t$-partitions with
$t = 1$.  
%% If all elements of $\mathcal{P}$ have dimensions in $\{s,s'\}$ with $s>s'$, one
%% also speaks of an $(s',s)$-spread or a partial $(s',s)$-spread, see \cite{metsch1999bose}.
Besides the simplicity %%genuineness 
of the proposed generalization, there are some similarities that promise interesting
applications.
The class of vector space $t$-partitions contains the $q$-analogs of a Steiner
systems, which are given by the cases where all elements of $\mathcal{P}$ have the
same dimension.
As a further generalization we mention the possibility of replacing {\lq\lq}contained in a
unique member of $\mathcal{P}${\rq\rq} by {\lq\lq}contained in exactly $\lambda$ elements of
$\mathcal{P}${\rq\rq}, which would include subspace designs ($q$-analogs of combinatorial
designs).
In the case $t = 1$, this has been considered in \cite{el2011lambda}.

Let $\mathcal{P}$ be a non-trivial vector space partition of $\F{v}{q}$ with a
non-empty set $\mathcal{N}$ of holes, and $k$ the 
second smallest dimension of the elements of $\mathcal{P}$. Then, we have
$\#\mathcal{N}\equiv \#\{N\in\mathcal{N}\,:\,N\le H\}\pmod {q^{k-1}}$ for each 
hyperplane $H$ of $\F{v}{q}$. This condition allows to conclude restrictions on
$\#\mathcal N$ independently of the dimension $v$ of the 
ambient space. Exploiting this congruence condition yields a series of improvements
\cite{kurz2016improved,nuastase2017maximum,kurz2017packing} 
for the maximum size of a partial $k$-spread, which is the set of the
$k$-dimensional elements of a vector space partition of type $k^{m_k} 1^{m_1}$. The 
underlying techniques can possibly be best explained using the language of
projective $q^{k-1}$-divisible linear codes and the linear programming 
method, see \cite{partial_spreads_and_vectorspace_partitions}.   

In our more general setting of a non-trivial vector space $t$-partition
$\mathcal{P}$ of $\F{v}{q}$, the set $\mathcal{N}$ of its $t$-dimensional subspaces
will play the role of the holes.
If $\mathcal{N} \neq\emptyset$ and $k>t$ is the second smallest dimension of the
elements of $\mathcal{P}$, we will prove $\#\mathcal{N}\equiv
\#\{N\in\mathcal{N}\,:\,N\le  H\}\pmod {q^{k-t}}$ for each hyperplane $H$ of
$\F{v}{q}$. Similarly,  we will study restrictions on $\#\mathcal{N}$ independently
of the dimension $v$ of the 
ambient space. Again some kind of linear programming method will be applied and
partially solved analytically.

The analog of partial $k$-spreads of maximum size are vector space $t$-partitions of
type $k^{m_k} t^{m_t}$ with maximum $m_k$ (for 
given parameters $v$ and $t$), which have been previously studied under the name
(optimal) constant-dimension codes. Denoting the maximum 
possible $m_k$ by $A_q(v,2k-2t+2;k)$, our main motivation for the introduction of 
vector space $t$-partitions are indeed the recent improvements of $257\le
A_2(8,6;4)\le 289$ to $257\le A_2(8,6;4)\le 272$ 
\cite{new_bounds_subspaces_codes} and finally $A_2(8,6;4)=257$ \cite{paper257}.
These parameters play a rather prominent role for constant-dimension codes, 
since, besides $A_2(6,4;3)=77<81$ \cite{hkk77}, all other known upper bounds for
$A_q(v,d;k)$, where $d\notin\{2k,2v-2k\}$, are obtained via the 
so-called Johnson bound and the existence of 
divisible codes, see \cite{kiermaier2017improvement}. While the result of
\cite{new_bounds_subspaces_codes} is based on more than 1000~hours computing time,
we 
will apply similar techniques in order to computationally show $m_3\le 240$ for all
vector space $2$-partitions of $\F{7}{2}$ of type 
$4^{17} 3^{m_3}2^{m_2}$ in less than seventy hours.
For a vector space $2$-partition of $\F{8}{2}$ of type $4^{m_4}2^{m_2}$ that
contains $17$ $4$-dimensional elements in a common hyperplane or passing through a
common point, this directly implies $m_4\le 257$, i.e., meeting the lower bound of
$A_2(8,6;4)$.
In the remaining cases, a direct counting argument gives $A_2(8,6;4)\le 272$.
 %% We also determine $m_3\le 247$ for vector space $2$-partions of $\F{7}{2}$ of type 
%% $4^{16}3^{m_3}2^{m_2}$, which implies $A_2(8,6;4)\le 263$.       

While the mentioned result is based on explicit computer computations for fixed
parameters, we have some hope that a more thorough study 
of vector space $2$-partitions may lead to an improvement of the currently best
known bound $A_q(8,6;4)\le \left(q^4+1\right)^2$ or of other parameters in general.
To that end 
we will present the first preliminary results on the existence of vector space
$t$-partitions and the possible cardinalities of the corresponding 
set $\mathcal{N}$ of $t$-dimensional subspaces satisfying $\#\mathcal{N}\equiv
\#\{N\in\mathcal{N}\,:\,N\le  H\}\pmod {q^{k-t}}$ for each 
hyperplane $H$ of $\F{v}{q}$. 

There is another connection with constant-dimension codes. Let $\mathcal{C}$ be the
set of $k$-dimensional elements of a vector space 
$(t+1)$-partition of $\F{2k}{q}$ of type $(k+t)^1k^\star (t+1)^\star$, where
$k>t+1\ge 1$. Replacing each element in $\mathcal{C}$ by its 
dual, we obtain a constant-dimension code in $\F{2k}{q}$ with minimum subspace
distance $2k-2t$ and cardinality $\#\mathcal{C}$ such that 
every codeword is disjoint from a $(k-t)$-dimensional subspace. 

Vector space $t$-partitions of type $k^{m_k} t^{m_t}$ are also of interest for the
set case, i.e., $q=1$. In other words, we are considering 
sets of $k$-subsets of $\{1,2,\dots,v\}$ such that every $t$-set is contained in
exactly one $k$-set (or contained in at most one $k$-set, 
if we anticipate the possible completion with $t$-sets). These structures are
equivalent to binary constant-weight codes with length 
$v$ and minimum Hamming distance $d\ge 2k-2t+2$. See e.g.\
\cite{agrell2000upper,best1978bounds} for upper bounds on $m_k$. 

The classification of the possible types of vector space $t$-partitions is also an
interesting problem for $q=1$. While it is trivial for 
$t=1$ it is not completely resolved for $t=2$. In the latter case one speaks of
pairwise balanced designs (with index $1$) or linear spaces, see e.g.\ 
\cite{colbourn2006handbook,dukes2015pairwise,ling1997pairwise}. In
\cite{spencer1968maximal} in has been shown that there is no set of triples covering
each pair 
exactly once except a single uncovered pair.\footnote{There exist e.g.\ $6$~triples
and $6$~quadruples of an $11$-set leaving exactly 
one pair uncovered and $12$ triples, $3$~quadruples, and a quintuple of a $12$-set
leaving exactly 
two intersecting pairs uncovered.} For more results in that direction we refer to
\cite{heinrich1992pairwise}.
%% Ich habe mir jetzt noch die beiden M??glichkeiten f??r zweielementige
%% Tails f??r q=1 t=2 angeschaut:
%% 
%% 1.
%% Zwei sich schneidende Paare (Vertreter {1,2}, {1,3}) geht erstmals f??r
%% v=12.
%% Die Konfiguration ist eindeutig, ihr Typ ist 2^2 3^12 4^3 5^1, ein
%% Vertreter ist
%% { 2, 3, 4, 5, 6 },
%% { 1, 4, 7, 10 }
%% { 1, 5, 8, 11 },
%% { 1, 6, 9, 12 },
%% { 2, 8, 9 },
%% { 2, 7, 11 },
%% { 2, 10, 12 },
%% { 3, 11, 12 },
%% { 3, 8, 10 },
%% { 3, 7, 9 },
%% { 4, 8, 12 },
%% { 5, 9, 10 },
%% { 6, 10, 11 },
%% { 4, 9, 11 },
%% { 5, 7, 12 },
%% { 6, 7, 8 },
%% { 1, 2 },
%% { 1, 3 }
%% 
%% 2.
%% Zwei disjunkte Paare (Vertreter {1,2}, {3,4}) geht erstmals f??r v=13.
%% Es gibt zwei Isomorphietypen:
%% 
%% 2a.
%% Typ 2^2 3^10 4^6 5^1
%% { 1, 3, 7, 8, 9 },
%% { 4, 6, 9, 11 },
%% { 1, 4, 5, 12 },
%% { 2, 3, 6, 13 },
%% { 2, 5, 9, 10 },
%% { 2, 7, 11, 12 },
%% { 4, 7, 10, 13 },
%% { 1, 6, 10 },
%% { 1, 11, 13 },
%% { 2, 4, 8 },
%% { 3, 5, 11 },
%% { 3, 10, 12 },
%% { 5, 6, 7 },
%% { 5, 8, 13 },
%% { 6, 8, 12 },
%% { 8, 10, 11 },
%% { 9, 12, 13 },
%% { 1, 2 },
%% { 3, 4 }
%% 
%% 2b.
%% Typ 2^2 3^14 4^4 5^1
%% { 5, 6, 7, 8, 9 },
%% { 1, 5, 10, 11 },
%% { 2, 7, 11, 13 },
%% { 3, 7, 10, 12 },
%% { 4, 5, 12, 13 },
%% { 1, 3, 8 },
%% { 1, 4, 7 },
%% { 1, 6, 13 },
%% { 1, 9, 12 },
%% { 2, 3, 5 },
%% { 2, 4, 8 },
%% { 2, 6, 12 }
%% { 2, 9, 10 },
%% { 3, 6, 11 },
%% { 3, 9, 13 },
%% { 4, 6, 10 },
%% { 4, 9, 11 },
%% { 8, 11, 12 },
%% { 8, 10, 13 },
%% { 1, 2 },
%% { 3, 4 }

The remaining part of this article is structured as follows. In
Section~\ref{sec_preliminaries} we introduce the preliminaries     
before we study the existence of vector space $t$-partitions in
Section~\ref{sec_existence}. As a contained substructure, $q^r$-divisible 
sets of $t$-subspaces are introduced and studied in Section~\ref{sec_divisible}. We
close with several open problems and a conclusion in 
Section~\ref{sec_conclusion}. 

\section{Preliminaries}
\label{sec_preliminaries}
\noindent
We briefly call $k$-dimensional subspaces of $\F{v}{q}$ \emph{$k$-subspaces}. $1$-subspaces are called points, $2$-subspaces are called lines, 
$3$-subspaces are called planes, $4$-subspaces are called solids, and $(v-1)$-subspaces are called hyperplanes. The number of $k$-subspaces 
in $\F{v}{q}$ is given by the Gaussian binomial coefficient $\gaussm{v}{k}{q}:=\prod_{i=0}^{k-1} \frac{q^{v-i}-1}{q^{k-i}-1}$. 

\begin{Definition}
  Let $t\in\mathbb{N}_{>0}$. A \emph{vector space $t$-partition of $\F{v}{q}$} is a set $\mathcal{P}$ of subspaces of $\F{v}{q}$ 
  such that every  $t$-subspace of $\F{v}{q}$ is contained in exactly one element of $\mathcal{P}$ and all elements of $\mathcal{P}$ have 
  dimension at least $t$ (so that they are incident with at least one $t$-subspace). We call $\mathcal{P}$ \emph{trivial} if all elements 
  either have dimension $t$ or $v$. If $\mathcal{P}$ contains $m_i$ elements of dimension $t\le i\le v$ we call 
  $v^{m_v} (v-1)^{m_{v-1}}\dots t^{m_t}$ the \emph{type of $\mathcal{P}$}, where $i^{m_i}$ can also be omitted if $m_i=0$.
\end{Definition}

As an example we consider vector space $2$-partitions of $\F{13}{2}$ of type $3^{1597245}$, which correspond to $2$-Steiner systems 
of planes in $\F{13}{2}$, whose existence has been proved in \cite{braun2016existence}. The existence of a vector space $2$-partition of $\F{7}{2}$ of type $3^{381}$ is 
equivalent to the existence of a binary $q$-analog of the Fano plane. If it exists it has an automorphism group of order at most 
two \cite{kiermaier2017order,braun2016automorphism}, the number of incidences between the blocks and other $k$-subspaces are 
known \cite{kiermaier2015intersection}, and not all sets of blocks incident with a point can correspond to a Desarguesian line 
spread \cite{heden2016existence,thomas1996designs}. Possible substructures 
of a $q$-analog of the Fano plane presently trigger a lot of research, see e.g.~\cite{braun2018q,hooker2017residual} and the references 
therein. The maximum known value 
of $m_3$ of a vector space $2$-partition of $\F{7}{2}$ of type $3^{m_3}2^{m_2}$ is $m_3=333$ \cite{paper333}. For general results on 
the existence of vector space $t$-partitions of $\F{v}{q}$ of type $s^{m_s}t^{m_t}$, also known as (partial) $(s,t)$-spreads, we refer the 
reader to e.g.~\cite{dentice2009bose,metsch1999bose}.

For two $k$-subspaces $U,W$ in $\F{v}{q}$ the \emph{subspace distance} is given by $d_S(U,W)=\dim(U+W)-\dim(U\cap W)=\dim(U)+\dim(W)-2\dim(U\cap W)
=2k-2\dim(U\cap W)$.
\begin{Definition}
  A constant-dimension code $\mathcal{C}$ of $\F{v}{q}$ of constant dimension $k$ and minimum subspace distance $d$ is a set of $k$-subspaces 
  such that the dimension of the intersection of any pair of $k$-subspaces is at most $\left\lfloor k-d/2\right\rfloor$. By $A_q(v,d;k)$ we 
  denote the corresponding maximum size, i.e., the number of $k$-subspaces. 
\end{Definition}

Each vector space $t$-partition $\mathcal{P}$ of $\F{v}{q}$ of type $k^{m_k}t^{m_t}$ is in $1$-to-$1$-correspondence to a constant-dimension code 
$\mathcal{C}=\{U\in\mathcal{P}\,:\, \dim(U)=k\}$ with minimum distance at least $2k-2t+2$, so that $m_k\le A_q(v,2k-2t+2;k)$. Note that 
by duality we have $A_q(v,d;k)=A_q(v,d;v-k)$. For known bounds, we 
refer to the online table \url{http://subspacecodes.uni-bayreuth.de} \cite{TableSubspacecodes}. As an example for constant-dimension codes we would 
like to mention lifted maximum rank distance (MRD) codes, see \cite{delsarte1978bilinear,gabidulin1985theory,roth1991maximum}.
%%  which were independently invented by Delsarte \cite{delsarte1978bilinear}, Gabidulin \cite{gabidulin1985theory}, and Roth \cite{roth1991maximum}.

\begin{Theorem}(see \cite{silva2008rank})
  For positive integers $k,d,v$ with $k\le v$, $d\le 2\min\{k,v-k\}$, and $d\equiv 0\pmod{2}$, the size of a lifted MRD code $\mathcal{C}$ of $k$-subspaces in $\F{v}{q}$ with 
  minimum distance at least $d$ is given by $$M(q,k,v,d):=q^{\max\{k,v-k\}\cdot(\min\{k,v-k\}-d/2+1)}.$$ Moreover, there exists a $(v-k)$-dimensional subspace 
  $U$ of $\F{v}{q}$ such that every element of $\mathcal{C}$ has trivial intersection with $U$. The set of $\left(\min\{k,v-k\}-d/2+1\right)$-subspaces that is 
  disjoint to $U$ is perfectly covered by the codewords.  
\end{Theorem}

\begin{Corollary}
  \label{cor_mrd_construction}
  For non-negative integers $k,t,v$ with $k\ge t+2$ and $v\ge 2k-t+1$, there exists a vector space $(t+1)$-partition of $\F{v}{q}$ of type 
  $(v-k+t)^1 k^m(t+1)^\star$, where $\log_q m=\max\{k,v-k\}\cdot(\min\{k,v-k\}-k+t+1)$.  
\end{Corollary}
\begin{Proof}
  Consider a lifted MRD code $\mathcal{C}$ of $k$-subspaces in $\F{v}{q}$ with minimum distance $d=2k-2t$. Let $U$ be the $(v-k)$-subspace that 
  has trivial intersection with the elements from $\mathcal{C}$.  
  Add a %%an arbitrary 
  $(v-k+t)$-subspace containing $U$, and complete the construction by adding uncovered $(t+1)$-subspaces. 
  %% Complete $\mathcal{C}$ by an arbitrary $(v-k+t)$-subspace containing $U$ to the 
  %% desired vector space $(t+1)$-partition. 
\end{Proof}
We remark that the construction also works for $v=2k-t$, where we obtain a vector space $(t+1)$-partition of $\F{v}{q}$ of type  
$k^{m+1}(t+1)^\star$ with $m=q^k$. 

\section{Existence of vector space $t$-partitions}
\label{sec_existence}
\noindent
In this section we will study the possible types of vector space $t$-partitions of $\F{v}{q}$ for small dimensions $v$. 
Here we will assume $t\ge 2$ and refer to the survey \cite{heden2012survey} for the case $t=1$. From that paper we also transfer 
the first conditions on the parameters $m_i$ of a vector space $t$-partition $\mathcal{P}$ of $\F{v}{q}$ of type $k^{m_k}\dots t^{m_t}$. 
Since every $t$-subspace is contained in a unique element in $\mathcal{P}$, we have
\begin{equation}
  \label{eq_c_1}
  \sum_{i=t}^k m_i\cdot\gaussm{i}{t}{q}=\gaussm{v}{t}{q},
\end{equation}
which is called \textit{packing condition} in \cite{heden2012survey} for $t=1$. This equation allows us to 
suppress the precise value of $m_t$ as done in Corollary~\ref{cor_mrd_construction}. Due to the dimension formula  
$\dim(U + V)=\dim(U)+\dim(V)-\dim(U\cap V)$, for any two subspaces $U$ and $V$ of $\F{v}{q}$, 
we have 
\begin{equation}
  \label{eq_c_2}
  m_i\le 1\,\text{ if }\,2i>v+t-1\quad\text{and}\quad m_im_j=0\,\text{ if }\,i+j>v+t-1
\end{equation}
for $t\le i<j\le v$. The specialization to $t=1$ is called \textit{dimension condition} in \cite{heden2012survey}. 

Known constructions of vector space $t$-partitions are given by lifted MRD codes. If $\mathcal{P}$ is a vector space $t$-partition 
of $\F{v}{q}$, $U$ an element of $\mathcal{P}$, and $\mathcal{P}'$ a vector space $t$-partition of $U$, then 
$\left(\mathcal{P}\backslash U\right)\cup \mathcal{P}'$ is also a vector space $t$-partition of $\F{v}{q}$. We call $\mathcal{P}'$ 
\textit{derived} from $\mathcal{P}$, matching the definition of a derived vector space partition for $t=1$.

From equations (\ref{eq_c_1}) and (\ref{eq_c_2}) we conclude that for $t\le v\le t+1$ each vector space $t$-partition $\mathcal{P}$ 
of $\F{v}{q}$ is trivial, i.e., either 
$\mathcal{P}=\{\F{v}{q}\}$ -- type $v^1$ -- or $\mathcal{P}$ is given by the $\gaussm{v}{t}{q}$ $t$-subspaces of $\F{v}{q}$. In the following we will 
consider the non-trivial vector space $t$-partitions only. For $v=t+2$ the dimension condition allows $m_{t+1}=1$ only, so that 
the packing condition gives type $(t+1)^1 t^{\gaussm{t+2}{t}{q}-\gaussm{t+1}{t}{q}}$. Here, $\mathcal{P}$ consists of an arbitrary 
$(t+1)$-subspace $U$ and all $t$-subspaces not contained in $U$. So far, all discussed cases are unique up to isomorphism. For $v=t+3$ we 
get $m_{t+2}\le 1$ and $m_{t+2}m_{t+1}=0$ so that we have type $(t+2)^1 t^{\gaussm{t+3}{t}{q}-\gaussm{t+2}{t}{q}}$ or  
type $(t+1)^{m_{t+1}}t^{m_t}$. In the latter case 
we have $m_{t+1}\le A_q(t+3,4;t+1)=A_q(t+3,4;2)$. The corresponding objects to $m_{t+1}=A_q(t+3,4;t+1)=A_q(t+3,4;2)$ are so-called 
(partial) line spreads of maximum size. If $t$ is odd, then $A_q(t+3,4;2)=\left(q^{t+3}-1\right)/(q^2-1)$, and 
$A_q(t+3,4;2)=\left(q^{t+3}-q^2(q-1)-1\right)/(q^2-1)$ otherwise, see e.g.~\cite{beutelspacher1975partial}. Here, there are several 
isomorphism types in general. So, using derived vector space $t$-partitions, in $\F{t+3}{q}$ there exist vector space $t$-partitions of type $(t+1)^i t^{\gaussm{t+3}{t}{q}-i\gaussm{t+1}{t}{q}}$ for 
all $0\le i\le A_q(t+3,4;2)$. For $v=t+4$ we conclude from the dimension condition that only the types $(t+3)^1 t^\star$, 
$(t+2)^1 (t+1)^{a} t^\star$, and $(t+1)^b t^\star$ might be possible for a non-trivial vector space $t$-partition. In the latter case we have 
$b\le A_q(t+4,4;t+1)=A_q(t+4,4;3)$. Since the current knowledge on $A_q(t+4,4;3)$ is rather limited, we mention the known bounds for $t=2$ only: 
$A_2(6,4;3)=77$ with precisely $5$ attaining isomorphism types and $q^6+2q^2+2q+1 \le A_q(6,4;3)\le \left(q^3+1\right)^2=q^6+2q^3+1$ for 
$q\ge 3$, see \cite{hkk77}. For type $(t+2)^1 (t+1)^{a} t^\star$ Corollary~\ref{cor_mrd_construction} gives a construction with $a=q^{2t+2}$, 
which is tight for $t=2$.

\begin{Lemma}
  \label{lemma_mrd_like_bound}
  If $\mathcal{P}$ is a vector space $2$-partition of $\F{v}{q}$ of type $(v-k+1)^1 k^a 2^\star$, where $k\ge 3$ and $v\ge 2k$, then 
  $a\le q^{2(v-k)}$.
\end{Lemma}    
\begin{Proof}
  Let $U$ be the unique $(v-k+1)$-subspace of $\mathcal{P}$. The number of lines disjoint from $U$ is given by 
  $$
    \gaussm{v}{2}{q}-\gaussm{v-k+1}{2}{q}-\frac{1}{q}\cdot\gaussm{v-k+1}{1}{q}\cdot\left(\gaussm{v}{1}{q}-\gaussm{v-k+1}{1}{q}\right)
    =q^{2(v-k)}\cdot\frac{q^{2k-1}-q^{k+1}-q^{k}+q^2}{(q^2-1)(q-1)}.
  $$
  Since each $k$-subspace $K$ of $\mathcal{P}$ intersects $U$ in exactly one point and the number of lines in $K$ disjoint from a 
  given point is given by
  $$
    \gaussm{k}{2}{q}-\gaussm{k-1}{1}{q}=\frac{q^{2k-1}-q^{k+1}-q^{k}+q^2}{(q^2-1)(q-1)},
  $$ 
  we have $a\le q^{2(v-k)}$. 
\end{Proof}
So the size of the construction from Corollary~\ref{cor_mrd_construction} is met for all cases where $t\in \{0,1\}$. %%\marginpar{$t>1$ ???} 
For $t=0$ the upper bound follows 
from counting the $k$-subspaces disjoint to a $(v-k)$-subspace. Removing the $(v-k)$-subspace gives precisely the lifted MRD codes with corresponding parameters. 
For $t=1$ the construction of Corollary~\ref{cor_mrd_construction} is far from being unique. We remark that the five isomorphism types 
of constant-dimension codes meeting the upper bound $A_2(6,4;3)=77$ each contain subsets of $64$ codewords that intersect a fixed solid in 
precisely a point. Moreover, there are exactly four isomorphism types of $64$ planes that intersect a fixed solid in 
precisely a point. The three ones that do not equal the lifted MRD code have automorphism groups of orders $24$, $16$, and 
$12$, respectively, and all can be extended to a constant-dimension code of cardinality $77$. %%\marginpar{Geometrische Struktur dieser Objekte?}
Lemma~\ref{lemma_mrd_like_bound} is also valid for the set case $q=1$, where it is tight.

%% \begin{Lemma}
%%   Let $\mathcal{P}$ be a vector space $2$-partition of $\F{6}{q}$ of type $4^1 3^a 2^\star$, then $a\le q^6+2q^3-q$
%% \end{Lemma}
%% \begin{Proof}
%%   The $\gaussm{4}{1}{q}=q^3+q^2+q+1$ points in the unique solid of $\mathcal{P}$ can be incident to at most $q^3$ planes of 
%%   $\mathcal{P}$. All other points can be incident to at most $q^3+1$ planes of $\mathcal{P}$. Since every plane contains $\gaussm{3}{1}{q}=q^2+q+1$ 
%%   points, the number of planes of $\mathcal{P}$ is at most 
%%   $$
%%     \left\lfloor\frac{\gaussm{4}{1}{q}\cdot q^3 + \left(\gaussm{6}{1}{q}-\gaussm{4}{1}{q}\right)\cdot\left(q^3+1\right)}{\gaussm{3}{1}{q}}\right\rfloor
%%     =\left\lfloor q^6+2q^3-q+\frac{q^2+q}{q^2+q+1}\right\rfloor=q^6+2q^3-q.
%%   $$
%% \end{Proof}
%% For $q=2$ this upper bound gives $a\le 78$, which can be improved to $a\le 76$ due to $A_2(6,4;3)=77$. Since none of the five isomorphism types of 
%% constant-dimension codes attaining cardinality $77$ allows the removal of a plane and the addition of a solid such that the solid intersects 
%% the planes in at most a point (\textbf{needs to be checked!}) we have $64\le a\le 75$ for the maximum value. 

For $v\ge t+5$ the situation gets rather involved, so that we assume $t=2$ and $v=7$ in the remaining part of this section. The dimension 
condition allows just the following types: $6^1 2^\star$, $5^1 3^{\tilde{m}_3}2^\star$, $4^{m_4} 3^{m_3} 2^{\star}$, and $3^{\bar{m}_3} 2^\star$, where 
$\bar{m}_3\le A_q(7,4;3)$ with e.g.\ $333\le A_2(7,4;3)\le 381$ and $6978\le A_3(7,4;3)\le 7651$, see \cite{TableSubspacecodes}. For the other 
parameterized cases we have $m_4\le A_q(7,6;4)=A_q(7,6;3)=q^4+1$ and $\tilde{m}_3\le q^8$, which is tight, see Corollary~\ref{cor_mrd_construction} 
and Lemma~\ref{lemma_mrd_like_bound}. Now, let us first look at constructions for the two maximal values for $m_4$.
\begin{Lemma}
  \label{lemma_construction_m_3_m4}
  For each prime power $q\ge 2$ there exist vector space $2$-partitions of $\F{7}{q}$ of type $4^{m_4} 3^{m_3} 2^{\star}$ 
  with $(m_4,m_3)=\left(q^4+1,q^{8}-q^4\right)$ and $(m_4,m_3)=\left(q^4,q^{8}-q^4+q^2+q+1\right)$. 
\end{Lemma}
\begin{Proof}
  Let $\mathcal{C}_8$ be a lifted MRD code of $q^8$ solids in $\F{8}{q}$ with minimum distance $6$ and $U$ be the unique solid having 
  trivial intersection to the elements from $\mathcal{C}_8$. For an arbitrary hyperplane $H$ of $\F{8}{q}$ that does not contain $U$ we 
  set $\mathcal{C}_7:=\left\{V\cap H\,:\,V\in \mathcal{C}_8\right\}$, so that $\mathcal{C}_7$ consists of $q^4$ solids and 
  $q^{8}-q^4$ planes. If $S$ is an arbitrary solid in $H$ that contains $U\cap H$, then $\mathcal{C}_7\cup S$ together with the uncovered 
  lines of $H$ gives a vector space $2$-partition of $H$ with type $4^{q^4+1}3^{q^{8}-q^4}2^\star$. For the other case, consider 
  $r=\gaussm{3}{1}{q}<\gaussm{4}{1}{q}$ arbitrary solids $S_1,\dots,S_{r}$ in $H$ containing $U\cap H$. Denoting 
  the $\gaussm{3}{2}{q}=r$ lines contained in $U\cap H$ by $L_1,\dots,L_r$, we choose $r$ planes $E_1,\dots,E_r$ such that $L_i\subseteq 
  E_i\subseteq S_i$. With this, $\mathcal{C}_7\cup\{E_i\,:\, 1\le i\le r\}$ can be completed by the uncovered lines to a vector space $2$-partition 
  of $H$ of type $4^{m_4} 3^{m_3} 2^{\star}$ with $(m_4,m_3)=\left(q^4,q^{8}-q^4+q^2+q+1\right)$.   
\end{Proof}

With respect to upper bounds for $m_3$ we consider the objects of $\mathcal{P}$ that are incident to a given point $P$. Modulo $P$ 
we obtain vector space partitions of $\F{6}{2}$ of type $3^{\bar{m}_3}2^{\bar{m}_2}1^\star$. The possible types have been completely 
classified, see e.g.~\cite{heden2012survey}. If $\bar{m}_3=3j+r$ with $j\in\mathbb{N}$ 
and $r\in\{0,1,2\}$, then $\bar{m}_2\le 21-5r+r^2-7j=:f(\bar{m}_3)$. 

\begin{Lemma}
  \label{lemma_ILP_computation}
  If $\mathcal{P}$ is a vector space $2$-partition of $\F{7}{2}$ of type $4^{m_4} 3^{m_3} 2^{\star}$, then 
  $m_3\le 240$ if $m_4=17$ and %%$m_3\le 247$ 
  $m_3\le 276$ if $m_4=16$.  
\end{Lemma}
\begin{Proof}
  Let $\mathcal{S}$ be a set of $16$ or $17$ solids in $\F{7}{2}$ pairwise intersecting in a point. By dualization we obtain 
  a set of $16$ or $17$ planes in $\F{7}{2}$ with trivial intersection. Those configurations have been classified up to symmetry in 
  \cite{honold2016classification}. Given all possible choices for $\mathcal{S}$, we develop an integer linear programming formulation 
  for the maximization of $m_3$. For each plane $E$ in $\F{7}{2}$ we introduce a variable $x_E\in\{0,1\}$ with $E\in\mathcal{P}$ iff 
  $x_E=1$, so that $m_3=\sum_{E\le\F{7}{2}} x_E$. If $L$ is a line of $\F{7}{2}$ that is contained in an element of $\mathcal{S}$, then we have 
  $\sum_{L\le E\le\F{7}{2}} x_E=0$ and $\sum_{L\le E\le\F{7}{2}} x_E\le 1$ otherwise. The LP relaxation of the current formulation can be further improved 
  by adding $\sum_{P\le E\le \F{7}{2}} x_E\le f(\tau(P))$, where $P$ is a point in $\F{7}{2}$ and $\tau(P)$ counts the number of elements of $\mathcal{S}$ 
  that contain $P$. Given $\mathcal{S}$ we denote the corresponding integer linear programming formulation by $ILP_{\mathcal{S}}$ and its 
  LP relaxation by $LP_{\mathcal{S}}$.
  
  For $\#\mathcal{S}=17$ it took 7~minutes to compute the $715$ linear programs $LP_{\mathcal{S}}$. Except $10$ cases, all of them have 
  a target value strictly less than $240$. In exactly one case a target value of $240$ can be attained for $ILP_{\mathcal{S}}$, which took 
  less than 66~hours to verify computationally. For $\#\mathcal{S}=16$ we computed the $14445$ instances $LP_{\mathcal{S}}$ leaving just $28$ 
  cases with a target value of at least $247$. 
  %% Es gibt 14.445 16er Konfigurationen. 14.417 haben einen Optimalwert < 247 in einer
  %% Stunde mit nur <=1 Bedingungen. Die übrigen 28 Indizes sind:
  %% 1,2,3,11,15,3557,3818,6026,6057,6111,6122,6130,6152,6184,6205,6216,6247,7498,7549,7560,7571,7672,7723,7764,7855,7876,7917,8793  
  It took 6~h to verify that $ILP_{\mathcal{S}}$ has a target value of at most $276$ %%$247$ 
  for these $28$~instances. After 99~h there remain just $7$~instances which may yield a target value strictly greater 
  than $247$, i.e., the lower bound given by Lemma~\ref{lemma_construction_m_3_m4}.
  %% Index        | Obergrenze Anzahl Planes nach max 99h
  %% 1            | 276
  %% 2            | 248
  %% 3            | 258
  %% 15           | 258
  %% 3557         | 262
  %% 7498         | 247
  %% 7876         | 251 
\end{Proof}

Let $a_i$ denote the number of points of $\F{7}{2}$ that are contained in exactly $i$ solids of $\mathcal{S}$. 
We remark that for $\#\mathcal{S}=17$, we can easily deduce $a_1=7$, $a_2=112$, and $a_3=8$, so that 
$m_3\le \frac{7\cdot f(1)+112\cdot f(2)+8\cdot f(3)}{7}=273$. For $\#\mathcal{S}<17$ even less information on the 
$a_i$ is sufficient to establish a competitive upper bound for $m_3$. 

\begin{Lemma}
  \label{lemma_upper_bound_m_4_m_3}
  If $\mathcal{P}$ is a vector space $2$-partition of $\F{7}{2}$ of type $4^{m_4} 3^{m_3} 2^{\star}$, then 
  $m_3\le 381-\left\lceil\frac{m_4(61-m_4)}{7}\right\rceil$.  
\end{Lemma}
\begin{Proof}
  %% Each $i$-subspace in $\mathcal{P}$ covers $\gaussm{i}{2}{q}$ lines. Since each line is covered at most once by the solids and planes 
  %% in $\mathcal{P}$, we have $m_4\gaussm{4}{2}{q}+m_3\gaussm{3}{2}{q}\le\gaussm{7}{2}{q}$.
  Let $a_i$ denote the number of points of $\F{7}{2}$ that are contained in exactly $i$ %%out 
  of the $m_4$ solids of $\mathcal{P}$. Counting points 
  gives $\sum_{i\ge 0} a_i=\gaussm{7}{1}{2}=127$ and $\sum_{i\ge 0} ia_i=\gaussm{4}{1}{2}m_4=15m_4$. Since every pair of solids of $\mathcal{P}$ 
  intersects in exactly one point, we additionally have $\sum_{i\ge 0} i(i-1) a_i=m_4(m_4-1)$. With this and the definition of the function $f$, 
  $\left\lfloor\frac{1}{7} \sum_{i\ge 0} f(i)\cdot a_i\right\rfloor$ is an upper bound for $m_3$. Next we maximize $\sum_{i\ge 0} f(i)\cdot a_i$ 
  for non-negative integers $a_i$ satisfying the three equations stated above. Since Lemma~\ref{lemma_ILP_computation} gives a 
  stronger bound than $m_3\le 274$ for $m_4=17$, we can assume $m_4\le 16$ in the following. From the last two equations we conclude
  $$
    a_1=m_4(16-m_4)+\sum_{i\ge 3} i(i-2)a_i\ge \sum_{i\ge 3} (2i-3)a_i,
  $$  
  so that $a_1\ge 2l-3$ if $a_l\ge 1$ for some $l\ge 3$. We claim that $a_i=0$ for all $i\ge 3$ in an optimal solution. 
  Assume $a_l\ge 1$ for some $l\ge 3$. Now, we modify the given $a_i$-vector by decreasing $a_l$ by $1$, increasing $a_{l-1}$ by $1$, increasing $a_2$ 
  by $l-1$, decreasing $a_1$ by $2l-3$ and increasing $a_0$ by $l-2$. The resulting vector $(a_0',a_1',\dots)$ satisfies the three equations 
  and has non-negative integer entries. By this operation the value of $\sum_{i\ge 0} f(i)\cdot a_i$ increases by $f(l-1)-f(l)+2l-6\ge f(l-1)-f(l)\ge 1$. 
  Thus, the optimal solution is given by $a_2={m_4\choose 2}$, $a_1=m_4(16-m_4)$, and $a_0=127-\frac{m_4(31-m_4)}{2}$ with 
  $$
    \left\lfloor \frac{1}{7}\sum_{i\ge 0} f(i)\cdot a_i\right\rfloor=
    \left\lfloor\frac{1}{7}\cdot\left(m_4^2-61m_4+2667\right)\right\rfloor=
    381-\left\lceil\frac{m_4(61-m_4)}{7}\right\rceil.
  $$
  %% For $m_4=17$, we automatically have $a_1=7$, $a_2=112$, $a_3=8$, and $a_i=0$ for all other $i$ due to the unique hole configuration 
  %% of a partial plane spread in $\F{7}{2}$ of cardinality $17$, which yields an upper bound of $273<274=
  %% 381-\left\lceil\frac{17(61-17)}{7}\right\rceil$. Later on, we give a computational upper bound of $240$ anyway.
\end{Proof}

%% \begin{eqnarray*}
%%   \left\lfloor\frac{15\cdot 17+(127-15)\cdot 21}{7}\right\rfloor &=& 372\\
%%   \left\lfloor\frac{1\cdot 15+2\cdot 14\cdot 17+(127-29)\cdot 21}{7}\right\rfloor &=& 364\\
%%   \left\lfloor\frac{1\cdot 14+3\cdot 14\cdot 17+(127-43)\cdot 21}{7}\right\rfloor &=& 356\\
%%   \left\lfloor\frac{3\cdot 15+3\cdot 13\cdot 17+(127-42)\cdot 21}{7}\right\rfloor &=& 356\\
%% \end{eqnarray*}

We remark that Lemma~\ref{lemma_upper_bound_m_4_m_3} gives $m_3\le 278$ for $m_4=16$. Summarizing the binary case $q=2$, we 
have the following bounds for $\max m_3$:
\begin{center}\tiny
  \begin{tabular}{rrrrrrrrrrrrrrrrrrr}
    \hline 
    $m_4$ & 17 & 16 & 15 & 14 & 13 & 12 & 11 & 10 & 9\\
    $\max m_3$ & 240 & 247\dots276 & 248\dots282 & 249\dots287 & 252\dots291 & 273\dots297 & 274\dots302 & 275\dots308 & 276\dots314\\
    \hline
    $m_4$ & 8 & 7 & 6 & 5 & 4 & 3 & 2 & 1 & 0\\    
    $\max m_3$ & 284\dots320 & 285\dots327 & 286\dots333 & 287\dots341 & 291\dots348 & 297\dots356 & 300\dots364 & 312\dots 372 & 333\dots 381\\
    \hline
  \end{tabular}
\end{center}
The upper bounds are obtained from Lemma~\ref{lemma_ILP_computation} and Lemma~\ref{lemma_upper_bound_m_4_m_3}. Lemma~\ref{lemma_construction_m_3_m4} 
gives constructions for $m_4\in\{16,17\}$. The construction for $m_4=0$ is taken from \cite{paper333}. For $m_4\in\{1,2,3,4,8,12,13\}$ the stated 
lower bounds are obtained from an integer linear programming formulation with prescribed subgroups of the automorphism group, i.e., the 
Kramer--Mesner approach. All other lower bounds are obtained by replacing a solid by a plane contained in the solid.

\section{$q^r$-divisible sets of $t$-subspaces}
\label{sec_divisible}
\noindent
Besides the conditions of Equation~(\ref{eq_c_1}) and Equation~(\ref{eq_c_2}), there is another technique for excluding the existence 
of (ordinary) vector space partitions, which just takes into account the second smallest occurring dimension. To this end, let  
$\mathcal{P}$ be a non-trivial vector space partition of $\F{v}{q}$, $\mathcal{N}\neq \emptyset$ be its set of holes\footnote{In, e.g., 
\cite{heden2009length} the author speaks of the tail of the vector space partition and considers lower bounds for its length, 
i.e., the cardinality of $\mathcal{N}$.}, i.e., $1$-dimensional 
elements, and $k$ be the second smallest dimension of the elements of $\mathcal{P}$. Then, we have 
$\#\mathcal{N}\equiv \#\{ N\in\mathcal{N}\,:\,N\le H\}\pmod {q^{k-1}}$ for each hyperplane $H$ of $\F{v}{q}$. Assigning a weight $w(H)\in\mathbb{N}$ 
to every hyperplane $H$ via $w_{\mathcal{N}}(H):=\#\mathcal{N}- \#\{ N\in\mathcal{N}\,:\,N\le H\}$, we can say that the weights of the hyperplanes 
are divisible by $q^{k-1}$. So, we also call the set $\mathcal{N}$ of points $q^{k-1}$-divisible. The possible cardinalities of 
$q^r$-divisible sets of points, or equivalently the length of $q^r$-divisible linear codes, see \cite{partial_spreads_and_vectorspace_partitions}, 
are quite restrictive. This approach allows to exclude the existence of vector space partitions without knowing the precise values 
of the $m_i$ or the dimension $v$ of the ambient space. The asymptotic result on the maximal cardinality of partial spreads from 
\cite{nuastase2017maximum} can e.g.\ be obtained using $q^r$-divisible sets of points, see \cite{kurz2017packing}. However, there are some 
rare cases where the existence of a vector space partition was excluded with more involved techniques, see e.g.\  \cite{el2010partitions} 
for the exclusion of a vector space partition of type $4^{13} 3^6 2^6$ in $\F{8}{2}$. Nevertheless, the classification of all possible 
cardinalities of $q^r$-divisible sets of points is an important relaxation. So far, in the binary case, the classification is complete 
for $r\le 2$ only, see \cite[Theorem 13]{partial_spreads_and_vectorspace_partitions}, while there is a single open case for $r=3$. 
A general result for \textit{small} cardinalities but arbitrary parameters $q$ and $r$ was obtained in 
\cite[Theorem 12]{partial_spreads_and_vectorspace_partitions}, see Theorem~\ref{thm_exclusion_q_r}. For each pair of parameters there is a 
largest integer $\operatorname{F}(q,r)$, 
called \emph{Frobenius number}, such that no $q^r$-divisible set of points of cardinality $\operatorname{F}(q,r)$ exists, see e.g.\ 
\cite{projective_divisible_binary_codes} for some bounds. For $q^r$-divisible multisets of points the possible cardinalities have been 
completely characterized in \cite{kiermaier2017improvement}.       

The aim of this section is to generalize the notion of $q^r$-divisible sets of points to $q^r$-divisible sets of $t$-subspaces and to deduce 
restrictions for the possible cardinalities of such sets.  

\begin{Definition}
  Let $\mathcal{C}$ %%\neq \emptyset$ 
  be a set of $t$-subspaces in $\F{v}{q}$. We call $\mathcal{C}$ $q^r$-divisible if 
  $\# \mathcal{C}\equiv \# \{C\in\mathcal{C}\,:\,C\le H\}\pmod {q^r}$ for all hyperplanes $H$ of $\F{v}{q}$. 
\end{Definition}

The link between $q^r$-divisible sets of $t$-subspaces and vector space $t$-partitions is given by: 

\begin{Proposition}
  Let $\mathcal{P}$ be a non-trivial vector space $t$-partition of $\F{v}{q}$ with $m_i=0$ for all $t<i<k$, %% and $m_t>0$, 
  then the set $\mathcal{N}$ of 
  $t$-subspaces of $\mathcal{P}$ is $q^{k-t}$-divisible.
\end{Proposition}
\begin{Proof}
  Using the convention $\gaussm{l-1}{0}{q}=1$, we have $\gaussm{l}{t}{q}-\gaussm{l-1}{t}{q}=\gaussm{l-1}{t-1}{q}\cdot q^{l-t}$, which is
  divisible by $q^{k-t}$ for all $l\ge k$. Note that we have $v> k$ since $\mathcal{P}$ is non-trivial. Counting the $t$-subspaces of $\F{v}{q}$ 
  gives $\sum_{i=k}^{v-1} m_i\gaussm{i}{t}{q}+\#\mathcal{N}=\gaussm{v}{t}{q}$. Now, let $H$ be an arbitrary hyperplane of $\F{v}{q}$, $\mathcal{N}'$ 
  be the set of elements of $\mathcal{N}$ that are contained in $H$, and $\mathcal{P}':=\{U\cap H\,:\, U\in\mathcal{P}, \dim(U)\ge k\}\cup\mathcal{N}'$ 
  be a vector space $t$-partition of $H$ of type $(v-1)^{m_{v-1}'}\dots (k-1)^{m_{k-1}'}(t)^{\#\mathcal{N}'}$, where we allow $t=k-1$, slightly 
  abusing notation. With this, we have $\sum_{i=k-1}^{v-1} m_i'\gaussm{i}{t}{q}+\#\mathcal{N}'=\gaussm{v-1}{t}{q}$. By subtracting both equations 
  we conclude that $q^{k-t}$ divides $\#\mathcal{N}-\#\mathcal{N}'$ since each $i$-subspace in $\mathcal{P}$ with $i\ge k$ corresponds either to 
  an $i$-subspace or an $(i-1)$-subspace in $\mathcal{P}'$ and $q^{k-t}$ divides $\gaussm{l}{t}{q}-\gaussm{l-1}{t}{q}$ for $l\ge k$.  
\end{Proof}

In the following let $\mathcal{N}$ be a $q^{r}$-divisible set of $t$-subspaces in $\F{v}{q}$ with minimal $v$, i.e., $\mathcal{N}$ is not 
completely contained in any hyperplane. By $a_i$ we denote the number of hyperplanes $H$ of 
$\F{v}{q}$ with $\#\{N\in\mathcal{N}\,:\,N\le H\}=i$ and set $n:=\#\mathcal{N}$. Double-counting the incidences of 
the tuples $(H)$ and $(B,H)$, where $H$ is a hyperplane and $B\in\mathcal{N}$ with $B\le H$ gives:
\begin{equation}
  \sum_{i=0}^{n-1} a_i=\gaussm{v}{1}{q} \quad\text{and}\quad
  \sum_{i=0}^{n-1} ia_i=n\cdot \gaussm{v-t}{1}{q}.\label{eq1}
\end{equation} 
For two different elements $B_1,B_2$ of $\mathcal{N}$ their span $\langle B_1,B_2\rangle$ has a dimension $i$ between $t+1$ and $2t$. Denoting 
the number of corresponding ordered pairs by $b_i$, double-counting gives:
\begin{equation}
  \sum_{i=0}^{n-1} i(i-1)a_i=\sum_{i=t+1}^{2t} b_i\gaussm{v-i}{1}{q} \quad\text{and}\quad \sum_{i=t+1}^{2t} b_i =n(n-1).\label{eq2}
\end{equation}

As a first non-existence criterion we state:

\begin{Lemma}
  \label{lemma_linear_condition}
  For a non-empty $q^r$-divisible set $\mathcal{N}$ of $t$-subspaces in $\F{v}{q}$, there exists a hyperplane $H$ with 
  $\#\{N\in\mathcal{N}\,:\,N\le H\}< n/q^t$, where $n=\#\mathcal{N}$.
\end{Lemma}
\begin{Proof}
  Let $i$ be the smallest index with $a_i\neq 0$. 
  Then, the equations of~(\ref{eq1}) are equivalent to $\sum_{j\ge 0} a_{i+q^rj}=\gaussm{v}{1}{q}$ and 
  $\sum_{j\ge 0} \left(i+q^rj\right)\cdot a_{i+q^rj}=n\gaussm{v-t}{1}{q}$. 
  Subtracting 
  $i$ times the first equation from the second equation gives $\sum_{j>0} q^rja_{i+q^rj}=n\cdot\frac{q^{v-t}-1}{q-1}-i\cdot\frac{q^v-1}{q-1}$. 
  Since the left-hand side is non-negative, we have $i\le \frac{q^{v-t}-1}{q^v-1}\cdot n< \frac{n}{q^t}$. 
\end{Proof}

The proof of Lemma~\ref{lemma_linear_condition} expresses the simple fact that a hyperplane with the minimum number of $t$-subspaces 
in $\mathcal{N}$ contains at most as many $t$-subspaces as the average number of $t$-subspaces in $\mathcal{N}$ per hyperplane. 
Lemma~\ref{lemma_linear_condition} excludes e.g.\ the existence of $q$-divisible sets $\mathcal{N}$ of $t$-subspaces in $\F{v}{q}$ of a 
cardinality $n\in[1,q-1]$. 

Next we turn to constructions of $q^r$-divisible sets of $t$-subspaces. 
For $t=1$ each $k$-subspace and each affine $k$-subspace, i.e., the difference of a $(k+1)$-subspace and a contained 
$k$-subspace, yields a $q^{k-1}$-divisible set. With this, the next lemma shows that $q^r$-divisible sets of $t$-subspaces of cardinality 
$q^{r+1}$ exist for all integers $r,t\ge 1$. 
 
\begin{Lemma}
  \label{lemma_construction_1}
  Let $\mathcal{N}$ be a $q^r$-divisible set of $t$-subspaces in $\F{v}{q}$ such that $q^r$ divides $\#\mathcal{N}$. Then, for each $s\in\mathbb{N}$ 
  there exists a $q^r$-divisible set $\mathcal{N}'$ of $(t+s)$-subspaces in $\F{v+s}{q}$. 
\end{Lemma} 
\begin{Proof}
  Assume $s\ge 1$, choose an $s$-subspace $U$ in $\F{v+s}{q}$ such that  $\F{v}{q}\oplus U=\F{v+s}{q}$, and set $\mathcal{N}'=\{U+N\,:\,N\in\mathcal{N}\}$.
\end{Proof}

\begin{Lemma}
  \label{lemma_construction_2}
  For integers $t\ge 1$ and $a\ge 2$ let $\mathcal{N}$ be a $t$-spread in $\F{at}{q}$, i.e., a set of $\frac{q^{at}-1}{q^t-1}$ disjoint $t$-subspaces. 
  Then $\mathcal{N}$ is $q^{(a-1)t}$-divisible.
\end{Lemma}
\begin{Proof}
  Since any point in $\F{at}{q}$ is contained in a unique member of $\mathcal{N}$ and $x\cdot \gaussm{t}{1}{q}+
  \left(\frac{q^{at}-1}{q^t-1}-x\right)\cdot\gaussm{t-1}{1}{q}=\gaussm{at-1}{1}{q}$ for $x=\frac{q^{(a-1)t}-1}{q^t-1}$, every hyperplane 
  contains exactly $x$ elements from $\mathcal{N}$. The divisibility follows from $\frac{q^{at}-1}{q^t-1}-\frac{q^{(a-1)t}-1}{q^t-1}=q^{(a-1)t}$.
\end{Proof}
We remark that $t$-spreads exist for all values of $t$, $a$, and $q$. Examples can e.g.\ be obtained from the so-called \textit{subfield construction}, 
i.e., taking all $\gaussm{a}{1}{q^t}=\frac{q^{at}-1}{q^t-1}$ points in $\F{a}{q^t}$ considering $\mathbb{F}_{q^t}$ as a $t$-dimensional vector space over $\mathbb{F}_q$.

\begin{Lemma}
  \label{lemma_construction_3}
  For integers $t\ge 1$, $s\ge 0$, and $a\ge 2$ let $\mathcal{N}$ be a union of $q^s$ disjoint $t$-spreads $\mathcal{S}_1,\dots,\mathcal{S}_{q^s}$ 
  in $\F{at}{q}$, i.e., $\mathcal{S}_i\cap\mathcal{S}_j=\emptyset$ for $i\neq j$. Then $\mathcal{N}$ is $q^{(a-1)t+s}$-divisible.
\end{Lemma}
\begin{Proof}
  For each hyperplane $H$ and each index $1\le i\le q^s$ we have $\# \mathcal{S}_i\equiv \#\{U\in S_i\,:\, U\le H\}\pmod {q^{(a-1)t}}$ 
  due to Lemma~\ref{lemma_construction_2}. The result follows from $\#\mathcal{N}=q^s\cdot\#\mathcal{S}_1$ and 
  $\#\{U\in \mathcal{N}\,:\, U\le H\}\equiv q^s\cdot \#\{U\in S_1\,:\, U\le H\}\pmod{q^{(a-1)t+s}}$.
\end{Proof}
In $\F{at}{q}$ there can be at most $\gaussm{at}{t}{q}\cdot\gaussm{t}{1}{q}/\gaussm{at}{1}{q}$ pairwise disjoint $t$-spreads, which is just the 
number of $t$-subspaces of $\F{at}{q}$ divided by the number of $t$-subspaces of a $t$-spread. If that upper bound is reached one speaks of a
\textit{$t$-parallelism}. These are known to exist for $(v=at,t,q)$ in $\{(2a,2,2),(2^i,2,q),(6,2,3),$ $(6,3,2)\}$ for all integers 
$a,i\ge 2$, see e.g.~\cite{etzion2015partial} and the cited references therein. So far, no non-existence result is known. If the stated upper 
bound on the number of $t$-spreads is not met, one speaks of a \textit{partial $t$-parallelism}. For the maximum number $p(v,t,q)$ of 
pairwise disjoint $t$-spreads in $\F{v}{q}$, the lower bounds $p(2a,2,q)\ge q^{2\left\lfloor \log (2a-1)\right\rfloor}+\dots+q+1$, 
$p(at,t,q)\ge 2^t-1$, and $p(at,t,q)\ge 2$, where $a\ge 2$, are proven in \cite{Beutelspacher1990} and \cite{etzion2015partial}. %%, respectively.   

Next we present a lower bound on the cardinality of a non-empty $q^r$-divisible set of $t$-subspaces:
\begin{Theorem}
  \label{thm_min_card}
  Let $t\ge 2$ and $r\ge 1$ %%$1\le r\le t$\textbf{???} 
  be integers and $\mathcal{N}\neq \emptyset$ be a $q^r$-divisible set of $t$-subspaces in $\F{v}{q}$, where $v$ is minimal. 
  \begin{enumerate}
    \item[(i)]  If $q^r$ divides $\#\mathcal{N}$, then $\#\mathcal{N}\ge q^{r+1}$. %%Moreover, if $\#\mathcal{N}= q^{r+1}$, then $\mathcal{N}$ 
                %%arises from the construction of Lemma~\ref{lemma_construction_1} applied to an affine $(r+1)$-subspace. %%\textbf{???}
    \item[(ii)] If $\#\mathcal{N}$ is not divisible by $q^r$, then $\#\mathcal{N}\ge q^{t}+1$ and $\#\mathcal{N}\ge q^r+ 
                \frac{q^{(\kappa-1)t}-1}{q^t-1}\cdot q^{r-(\kappa-1)t}$, where $\kappa$ is the smallest positive integer satisfying 
                $\frac{q^{\kappa t}-1}{q^t-1}\ge q^r$.                 
  \end{enumerate}
\end{Theorem}
\begin{Proof}
  \begin{enumerate}
    \item[(i)]
  Assume $\#\mathcal{N}=lq^r$ for some positive integer $l$. Setting $\Delta=q^r$, $y=q^{v-t-1}$, and $c_i=a_i(q-1)$ for all $0\le i\le \#\mathcal{N}-1$, 
  the equations from (\ref{eq1}) are equivalent to 
  $$
    \sum_{i=0}^{l-1} c_{i\Delta}=q^{t+1}y-1\quad\text{and}\quad
    \sum_{i=0}^{l-1} i(\Delta-1)c_{i\Delta}=l(\Delta-1) \left(qy-1\right).
  $$  
  From Equation~(\ref{eq2}) we conclude
  $$
    l(l\Delta-1)(q^{-t+1}y-1)\le \sum_{i=0}^{l-1} i(i\Delta-1)c_{i\Delta}\le l(l\Delta-1)(y-1), 
  $$  
  so that
  $$
    l(\Delta-1) \left(qy-1\right)=\sum_{i=0}^{l-1} i(\Delta-1)c_{i\Delta}\le \sum_{i=0}^{l-1} i(i\Delta-1)c_{i\Delta}\le l(l\Delta-1)(y-1).
  $$
  Since $l\ge 1$, we have $(\Delta-1) \left(qy-1\right)\le (l\Delta-1)(y-1)$, so that $q\Delta+\Delta y+y\le 2\Delta+qy$ for $l\le q-1$. Since $q\ge 2$ 
  and $\Delta\ge q$, we obtain $y\le 0$, which is a contradiction. Thus, $l\ge q$ and $\#\mathcal{N}\ge q^{r+1}$. 
  %% For the special case $n=q^{r+1}$ we have 
  %% $l=q$ and $(\Delta-1) \left(qy-1\right)\le (l\Delta-1)(y-1)$ implies $y\ge \Delta$. \textbf{\dots}
  \item[(ii)] 
  Assume $\#\mathcal{N}=lq^r+x$ with $0<x<q^r$ for some integers $x,l$. Lemma~\ref{lemma_linear_condition} gives $\#\mathcal{N}\ge q^t+1$ 
  and from the divisibility we conclude $l\ge 1$, so that we assume $l=1$ in the following. With this, $\Delta=q^r$, and $y=q^{v-t}$, the 
  equations from (\ref{eq1}) are equivalent to  
  $$
    x(q-1)a_x=x\left(q^ty-1\right)\quad\text{and}\quad
    x(q-1)a_x=(x+\Delta)(y-1),
  $$
  so that $\Delta/y=x+\Delta-xq^t$ and $0\le v-t\le r$. Isolating $x$ gives 
  $\left(q^t-1\right)x=(y-1)\cdot\frac{\Delta}{y}=\Delta\cdot\left(1-\frac{1}{y}\right)$, which implies that $q^t-1$ divides $y-1$, i.e., 
  $t$ divides $v$, and that $x$ is increasing with $y$. So, let $v=\kappa\cdot t$ for some positive integer $\kappa$ with $(\kappa-1)t\le r$. Then, 
  $x=\frac{q^{(\kappa-1)t}-1}{q^t-1}\cdot q^{r-(\kappa-1)t}$ is increasing with $\kappa$. Of course 
  $\#\mathcal{N}\le \gaussm{v}{1}{q}/\gaussm{t}{1}{q}$, so that $q^r\le \frac{q^{\kappa t}-1}{q^t-1}$.  
  
  \end{enumerate}  
\end{Proof}

The construction of Lemma~\ref{lemma_construction_1} and the remark before Lemma~\ref{lemma_construction_1} show that (i) is tight. If $r\le t$, then the first part of (ii) is tight due to the 
construction of Lemma~\ref{lemma_construction_2} with $a=2$. If $t$-parallelisms exist for all parameters (the dimension $v$ has of course 
to be divisible by $t$), then also the second part of (ii) is tight. The construction of Lemma~\ref{lemma_construction_3} shows that also 
a weaker assumption suffices for this claim. 

We remark that Theorem~\ref{thm_min_card} generalizes a theorem on the so-called length of the tail of a vector space partition, originating 
from \cite{heden1986partitions}, for the special case $t=1$, where the $k$-subspaces automatically are disjoint.

\begin{Theorem} (\cite[Theorem 10]{kurz2017heden})
  \label{thm_tail}
  For a non-empty $q^r$-divisible set $\mathcal{N}$ of pairwise disjoint $k$-subspaces in $\F{v}{q}$
  the following bounds on $n=\#\mathcal{N}$ are tight.
  \begin{enumerate}[(i)]
    \item We have $n\ge q^k+1$ and if $r\ge k$ then either $k$ divides $r$ and $n\ge \frac{q^{k+r}-1}{q^k-1}$ or $n\ge\frac{q^{(a+2)k}-1}{q^k-1}$, 
          where $r=ak+b$ with $0<b<k$ and $a,b\in\mathbb{N}$.
    \item Let $q^r$ divide $n$. If $r<k$ then $n\ge q^{k+r}-q^k+q^r$ and $n\ge q^{k+r}$ otherwise. 
  \end{enumerate}  
  %%All bounds are tight.
\end{Theorem} 

Aiming at characterizations of all possible cardinalities of $q^r$-divisible sets of $t$-subspaces it is useful to collect some 
more constructions. Taking the set of all $t$-subspaces gives another construction of divisible sets of $t$-subspaces.
\begin{Lemma}
  \label{lemma_construction_4}
  For integers $t\ge 1$ and $v\ge t+1$ the set $\mathcal{N}$ of all $t$-subspaces of $\F{v}{q}$ is $q^{v-t}$-divisible.
\end{Lemma} 
\begin{Proof}
  We have $\gaussm{v}{t}{q}-\gaussm{v-1}{t}{q}=\gaussm{v-1}{t-1}{q}\cdot q^{v-t}$.
\end{Proof}

%% Embedding two $q^r$-divisible sets of $t$-subspaces in orthogonal parts of a suitable ambient space, we obtain:
The set of achievable cardinalities of $q^r$-divisible sets of $t$-subspaces is closed under addition:
\begin{Lemma}
  \label{lemma_addition}
  Let $\mathcal{N}_1$ and $\mathcal{N}_2$ be $q^r$-divisible sets of $t$-subspaces in $\F{v_1}{q}$ and $\F{v_2}{q}$, respectively. 
  Then, there exists a $q^r$-divisible set of $t$-subspaces in $\F{v_1}{q}\oplus\F{v_2}{q}\cong\F{v_1+v_2}{q}$ with cardinality 
  $\#\mathcal{N}_1+\#\mathcal{N}_2$. 
\end{Lemma} 
In many cases an ambient space of dimension smaller than $v_1+v_2$ is sufficient.

%% For $q=2$, $r=1$, and $t=2$, the previous 
%% results allow to completely characterize the set $S_{r,q}^t$ of possible cardinalities of $q^r$-divisible sets of $t$-subspaces. From 
%% Theorem~\ref{thm_min_card} we conclude $S_{1,2}^2\subseteq \mathbb{N}_{\ge 4}$. Lemma~\ref{lemma_construction_1} gives examples 
%% for cardinalities $4$ and $6$, based on an affine plane and two disjoint lines. Lemma~\ref{lemma_construction_2} gives an example 
%% for cardinality $5$, based on a line spread of $\F{4}{2}$. The set of all seven lines in $\F{3}{2}$ yield an example for 
%% cardinality $7$, see Lemma~\ref{lemma_construction_4}, so that Lemma~\ref{lemma_addition} gives $S_{1,2}^2= \mathbb{N}_{\ge 4}$.

\section{Conclusion and open problems}
\label{sec_conclusion}
\noindent
Vector space $t$-partitions have many properties in common with ordinary vector space partitions, so that this class forms 
an interesting generalization. We have presented a few initial results on the existence of vector space $t$-partitions and their 
relaxation to $q^r$-divisible sets of $t$-subspaces. Only scratching the surface in this paper, we close with some open problems. 

While Lemma~\ref{lemma_mrd_like_bound} gives an upper bound on the cardinality of constant-dimension codes of dimension $k$ in $\F{2k}{q}$ 
with subspace distance $2k-2$ such that the codewords are disjoint from a $(k+1)$-subspace $U$, the underlying question is more general. 
What about $t>1$ in Lemma~\ref{lemma_mrd_like_bound}? If we forgo the link to vector space $t$-partitions via duality, we can ask for an 
upper bound on the cardinality of constant-dimension codes of dimension $k$ in $\F{v}{q}$ 
with subspace distance $d$ such that the codewords are disjoint from an $s$-subspace $U$. For the parameters $q=2$, $v=7$, $k=3$, $d=4$, and $s=3$ 
the corresponding lifted MRD code gives an example of cardinality $256$. So far we are only able to prove an upper bound of $290$.\footnote{Since no vector 
space partition of $\F{6}{2}$ of type $3^12^{18}1^2$ exists, every point $P$ (outside of $U$) can be contained in at most $17$ planes, which implies 
an upper bound of $\left\lfloor (127-7)\cdot 17/7\right\rfloor=291$. This upper bound can not be attained, since otherwise the argument from 
\cite{kiermaier2017improvement} gives a $4$-divisible multiset of $3$ points, which does not exist.} So, we ask for tighter bounds in this specific 
case and for the general problem.

In Section~\ref{sec_divisible} we have seen that the set of holes of a vector
space $t$-partition has to be a $q^r$-divisible set of $t$-subspaces. This
significantly restricts the possible types of vector space $t$-partitions
and raises the question how tight the resulting restrictions are.
For $q=1$, the condition of $q^r$-divisibility is trivially satisfied in
all cases. Indeed, we are not aware of any example of a hole
configuration $\mathcal{N}$ of $t$-subsets which provably is not realizable as a
vector space $t$-partition for $q=1$, i.e., a partition of the set of
$t$-subsets of a set $V$ such that all parts of size $t$ are given by $\mathcal{N}$.

%% That might be seen as an indication that the concept of
%% $q^r$-divisibility provides a reasonable condition on the realizability
%% of a given hole configuration as a $t$-vector space partition for general $q$.

Having determined the minimum possible cardinality of a $q^r$-divisible set of $t$-subspaces, for many parameters with $t\ge 2$, in 
Theorem~\ref{thm_min_card}, one can ask for the spectrum of possible cardinalities. For $t=1$ the following is known:
\begin{Theorem}(\cite[Theorem 12]{partial_spreads_and_vectorspace_partitions})
  \label{thm_exclusion_q_r}
  For the cardinality $n$ of a $q^r$-divisible set $\mathcal{C}$ of $1$-subspaces over
  $\mathbb{F}_q$ we have
  $$
    n\notin\left[(a(q-1)+b)\gaussm{r+1}{1}{q}+a+1,(a(q-1)+b+1)\gaussm{r+1}{1}{q}-1\right],
  $$
  where $a,b\in\mathbb{N}_0$ with $b\le q-2$, $a\le r-1$, and $r\in\mathbb{N}_{>0}$.
  
  In other words, if $n\le rq^{r+1}$, 
  then $n$ can be written as
$a\gaussm{r+1}{1}{q}+bq^{r+1}$ for some $a,b\in\mathbb{N}_0$.
\end{Theorem}

For $q=2$, $t=2$, and $r=1$ we remark that the possible cardinalities are given by $\mathbb{N}_{\ge 4}$. Examples of cardinality $4$ and 
$6$ are given by Lemma~\ref{lemma_construction_1}, Lemma~\ref{lemma_construction_2} gives a construction for cardinality $5$, and 
Lemma~\ref{lemma_construction_4} gives a construction for cardinality $7$, so that Lemma~\ref{lemma_addition} continues these constructions 
to all cardinalities in $\mathbb{N}_{\ge 8}$. For other parameters there are gaps in the sets of possible cardinalities. For which parameters 
can these sets be completely determined? What is the second smallest cardinality? Can Theorem~\ref{thm_exclusion_q_r} be generalized, i.e., 
for which ranges do integer combinations of two base constructions explain all possible cardinalities? What is the largest cardinality $n$ such 
that no $q^r$-divisible set of $t$-subspaces of cardinality $n$ exists? This number was called \textit{Frobenius number} for the special case $t=1$ 
in \cite{partial_spreads_and_vectorspace_partitions}. Determine bounds on the Frobenius number.  

Almost the same questions can be asked for vector space $t$-partitions. As for ordinary vector space partitions, the classification of 
all possible types is indeed a very hard problem if the dimension is not too small. However, for vector space $t$-partitions in $\F{7}{2}$ 
some improvements of the presented results seem to be achievable. Triggered by the motivating example of $A_2(8,6;4)<289$, we ask for a 
computer-free proof of $A_q(8,6;4)<\left(q^4+1\right)^2$. Nevertheless having just a very tiny numerical evidence, we state the following 
two rather strong conjectures in order to stimulate the search for counter examples.

\begin{Conjecture}
  $A_q(2k,2k-2;k)=q^{2k}+1$ for each $k\ge 4$.
\end{Conjecture}

We remark that the conjecture is true for the set case $q=1$, while $A_1(6,4;3)=2=1^6+1$ (slightly abusing notation).  

\begin{Conjecture}
  If $\mathcal{P}$ is a vector space $2$-partition of $\F{2k-1}{q}$ of type $k^{q^k+1} (k-1)^{m_{k-1}} 2^{\star}$, then 
  $m_{k-1}\le q^{2k}-q^k$ for all $k\ge 4$.
\end{Conjecture}

Again the conjecture is true for the set case $q=1$.

%%\bibliographystyle{amsplain}
%%\bibliographystyle{abbrv}
%%\bibliography{generalized_vector_space_partitions}

\end{document}